\documentclass[reqno]{amsart}
\theoremstyle{plain}
\newtheorem{lemma}{Lemma}
\newtheorem{theorem}[lemma]{Theorem}

\theoremstyle{remark}
\newtheorem{remark}{Remark}

\def\ee{\epsilon}
\def\aa{\alpha}

\def\om{\omega}

\def\tt{\theta}

\def\la{\lambda}
\def\Dd{\Delta_h}
\def\Gg{\Gamma}
\def\Om{\Omega}
\def\om{\omega}
\def\pp{\partial}

\begin{document}

\title[NS Equations]
{Pressure regularity criterion for the three
dimensional Navier--Stokes equations in infinite channel}

\date{February 22, 2007}
\thanks{\textit{ }}

\author[C. Cao]{Chongsheng Cao}
\address[C. Cao]
{Department of Mathematics  \\
Florida International University  \\
University Park  \\
Miami, FL 33199, USA}
\email{caoc@fiu.edu}

\author[E.S. Titi]{Edriss S. Titi}
\address[E.S. Titi]
{Department of Mathematics \\
and  Department of Mechanical and  Aerospace Engineering \\
University of California \\
Irvine, CA  92697-3875, USA \\
{\bf also}  \\
Department of Computer Science and Applied Mathematics \\
Weizmann Institute of Science  \\
Rehovot 76100, Israel} \email{etiti@math.uci.edu  and edriss.titi@weizmann.ac.il}

\begin{abstract}
In this paper we consider the three--dimensional  Navier--Stokes
equations in an infinite channel. We provide a sufficient condition,
in terms of $\partial_z p$, where $p$ is the pressure, for the
global existence of the strong solutions to the three--dimensional
Navier--Stokes equations.
\end{abstract}

\maketitle



\section{Introduction}   \label{S-1}

The question of global regularity for the $3D$ Navier--Stokes
equations is a major open problem in applied analysis. Over the
years there has been an intensive work by many authors attacking
this problem (see, e.g.,  \cite{CP95}, \cite{CF88}, \cite{DJ95},
\cite{LADY69}, \cite{LADY03}, \cite{LR02}, \cite{LL69}, \cite{PL96},
\cite{SO01A}, \cite{TT84},  \cite{TT95}, \cite{TT00} and references
therein). It is well known that the 2D Navier--Stokes equations have
a unique weak and strong solutions which exist globally in time
(cf., for example, \cite{CF88}, \cite{TT84}, \cite{TT00}). In the
$3D$ case, the weak solutions are known to exist globally in time.
But, the uniqueness, regularity, and continuous dependence on
initial data for weak solutions are still open problems.
Furthermore, strong solutions in the $3D$ case are known to exist
for a short interval of time whose length depends on the physical
data of the initial--boundary value problem. Moreover this strong
solution is known to be unique (cf., for example, \cite{CF88},
\cite{SO01A}, \cite{TT84}). The subtle difference between the $2D$
and $3D$ incompressible Navier--Stokes equations manifests itself in
a clear way in the vorticity formulation of these equations. Let $u$
denote the velocity field and $\om = \nabla \times u$ the vorticity.
The equations that govern the evolution of the vorticity are given
by
\begin{eqnarray*}
&&\hskip-.8in
\frac{\pp \om}{\pp t} - \nu  \Delta \om  + (u\cdot \nabla) \om -
\om \cdot \nabla u = 0,    \\
&&\hskip-.8in
 \nabla \cdot u =0.
\end{eqnarray*}
In this formulation the main obstacle for  proving the global
regularity is the vorticity stretching term $\om \cdot \nabla u$.
This term is identically equal to zero in the $2D$ case.
Nonetheless, there are some results regarding the global regularity
for the $3D$ Navier--Stokes equations under special symmetries and
for which the vorticity stretching term is non--trivial. For
example, the case of $3D$ axi--symmetric flows in domains of
revolution (with a positive distance from the $z-$axis)
\cite{LADY69} and \cite{LADY70} and helical flows \cite{LMT90}. It
is worth mentioning, however, that the questions of the global
regularity for $3D$ Euler equations (i.e., when the viscosity
$\nu=0$) for axi--symmetric and helical flows with nontrivial
vorticity stretching term are still open.

In addition, there are few known results regarding the global
regularity for special types of initial data. For instance, for
small $H^1$ initial data (cf. e.g., \cite{CF88}, \cite{TT84}).
Fujita and Kato \cite{HK64} proved the global well-posedness  for
small $H^{1/2}$ initial data, and later Kato \cite{KT84} proved the
same result for  small $L^3$ initial data  (see also \cite{GI86}).
An interesting global existence result for large, but very
``oscillatory", initial data was also proved in \cite{BV95},
\cite{BV96}, \cite{CM95}, \cite{MY06} (see also references therein).
The latter case can be regarded, roughly speaking, as small initial
data in $H^{1/2}$ and hence falls as a special case of results in
\cite{HK64}. The proofs of the above result relies heavily on the
viscosity mechanism. Roughly speaking, the viscosity acts very
strongly in dissipating the solution starting from  such oscillatory
initial data so that in a very short time the solution becomes very
small. As a result from that moment on a small initial data argument
applies to establish the global existence.  Here again, it is worth
stressing that nothing is known about global regularity for the $3D$
Euler equations for this kind of oscillatory initial data. It is
also worth mentioning the result of \cite{RS93} about global
existence of the strong solutions to the three-dimensional
Navier--Stokes equations for large initial data in thin domains,
where the upper bound on the size of the initial data depends
inversely on the thickness of the domain.

Taking a physical point view, we define the dimensionless Reynolds
number $Re= \frac{\nu U_0 }{\ell}$, where $U_0$ represents a typical
magnitude of velocity (e.g. the size of the initial data) and $\ell$
a typical length scale of the domain (e.g. its thickness). As a
result one can view, roughly speaking, the above mentioned global
regularity results to hold under the assumption of small enough
Reynolds number.

The other direction concerning the question of global regularity for
the $3D$ Navier--Stokes equations is to provide sufficient
conditions for the  global regularity. For example we refer the
reader to the pioneer work of Prodi \cite{PG59} and of Serrin
\cite{SJ62} (see also the survey paper of Ladyzhenskaya
\cite{LADY03} and references therein). Most recently, there has been
some progress along these lines which states that a strong solution
$u$ exists on the time interval $[0, T]$ as long as
\begin{equation}
u \in L^{\aa} ([0, T], L^{\beta} (\Om)), \quad\text{with}\quad
\frac{2}{\aa}+\frac{3}{\beta}=1, \quad \text{for} \quad \beta > 3,
\label{SP}
\end{equation}
(see, for example, \cite{BL02}, \cite{BG02}, \cite{GI86}, \cite{GM85},
\cite{KT84}, \cite{SO01}, \cite{SO02}, and references therein). Moreover, there are some sufficient
regularity conditions only on one component of the velocity field
of the $3D$ NSE on the whole space $\mathbb{R}^3$ or under
periodic boundary conditions (cf. e.g., \cite{HE02}, \cite{KZ05},
\cite{PO03}, \cite{ZY02}). In \cite{CQT05} we introduced a
sufficient regularity condition on one component of the velocity
field of the $3D$ Navier--Stokes equations under Dirichlet
boundary conditions (see also, \cite{MI06}).

\vskip0.1in

The Navier--Stokes equations of viscous incompressible fluid in an
infinite channel $\mathbb{R}^2 \times (0, 1) \subset \mathbb{R}^3$
read:
\begin{eqnarray}
&&\hskip-.8in
\frac{\pp v}{\pp t} - \nu  \Dd v -\nu v_{zz} + (v\cdot \nabla_h) v +w v_z
+ \nabla_h p = f,  \label{EQ-1}    \\
&&\hskip-.8in
\frac{\pp w}{\pp t} - \nu  \Dd w-\nu w_{zz} + v \cdot \nabla_h w +ww_z+ p_z = g,
\label{EQ-2}    \\
&&\hskip-.8in \nabla_h \cdot v +w_z =0,   \label{EQ-3}
\end{eqnarray}
where $v=(v_1, v_2)$, the horizontal velocity field components, $w$,
the vertical velocity component, and $p$, the pressure are unknowns.
$(f, g)$, the forcing term, and  $\nu > 0$, the viscosity are given.
We set $\nabla_h  = (\pp_x, \pp_y)$ to be the horizontal gradient
operator and $\Dd = \pp_x^2 +\pp_y^2$ is the horizontal Laplacian.
We denote by:
\begin{eqnarray}
&&\hskip-.8in
\Gg_u = \{ (x,y,1) \in \mathbb{R}^3 \},  \\
&&\hskip-.8in
\Gg_b = \{ (x,y,0) \in \mathbb{R}^3 \},
\end{eqnarray}
the physical solid boundaries of the channel $\mathbb{R}^2 \times (0, 1)$.
We equip the system (\ref{EQ-1})--(\ref{EQ-3}) with the following
no--normal flow and stress--free boundary conditions on the physical
boundaries of the channel, namely,
\begin{eqnarray}
&&\hskip-.8in  \; \frac{\pp v }{\pp z} =0, \; w=0;  \qquad \mbox{on }
\Gg_u \quad \mbox{and  } \; \; \Gg_b.
\label{B-1}
\end{eqnarray}
Furthermore, we assume that $(v, w)$ is periodic with period $1$
in both horizontal directions. That is,
\begin{equation}
v(x+1,y,z)=v(x,y+1,z)=v(x,y,z),   \quad
w(x+1,y,z)=w(x,y+1,z)=w(x,y,z). \label{B-2}
\end{equation}
Because of the horizontal periodic boundary conditions we consider
here \[
\Om=(0, 1)^3
\]
as our basic domain of the flow, and we denote
by \[
M=(0, 1)^2.
\]
In addition, we supply the system with the
initial condition:
\begin{eqnarray}
&&\hskip-.8in v(x,y,z,0) = v_0 (x,y,z),  \label{INIT-1}  \\
&&\hskip-.8in w(x,y,z,0) = w_0 (x,y,z).  \label{INIT-2}
\end{eqnarray}

Let us denote by $L^q(\Om)$ and $H^m (\Om)$ the usual
$L^q-$Lebesgue and Sobolev spaces, respectively (cf. \cite{AR75}). We
denote by
\begin{equation}
\| \phi\|_q =  \left(  \int_{\Om} |\phi|^q \; dxdydz
\right)^{\frac{1}{q}},   \qquad  \mbox{ for every $\phi \in
L^q(\Om)$}.
 \label{LQ}
\end{equation}
Let
\begin{eqnarray*}
 \mathcal{V}_1 &=&  \left\{ \phi \in \left( C^{\infty}(\mathbb{R}^2 \times
[0, 1]) \right)^2: \left. \frac{\pp \phi}{\pp z } \right|_{z=0}= 0; \left.
\frac{\pp \phi}{\pp z} \right|_{z= 1}= 0;
\phi(x+1,y,z)=\phi(x,y+1,z)=\phi(x,y,z) \right\},   \\
\mathcal{V}_2 &=&  \left\{ \psi \in C^{\infty}(\mathbb{R}^2 \times
[0, 1]): \left. \psi \right|_{z=0}= 0; \left.\psi \right|_{z= 1}=
0; \psi(x+1,y,z)=\psi(x,y+1,z)=\psi(x,y,z) \right\},  \\
 \mathcal{V} &=&  \left\{ (\phi, \psi): \phi \in  \mathcal{V}_1, \; \psi \in
 \mathcal{V}_2, \;
 \nabla_h \cdot \phi + \psi_z = 0 \right\}.
\end{eqnarray*}
We denote by $H$ and $V$ be the closure spaces of $\mathcal{V}$ in
$L^2(\Om)$ under the $L^2-$topology, and  in $H^1(\Om)$ under the
$H^1-$topology, respectively. We say $(v, w)$ is a Leray--Hopf weak
solution to the system (\ref{EQ-1})--(\ref{INIT-2}) if $(v, w)$
satisfies

\begin{itemize}

\item[(1)]
$(v, w)  \in C([0, T], H) \cap L^2([0, T], V),  $ and
$ (\pp_t v, \pp_t w)  \in L^1([0, T], V^{\prime}),$
where $V^{\prime}$ is the dual space of $V$;

\item[(2)] the weak formulation:
\begin{eqnarray*}
&&\hskip-0.35in
\int_{\Om} \left( v \cdot \phi + w \, \psi \right)\, dxdydz
-\int_{\Om} \left( v(t_0) \cdot \phi(t_0) + w(t_0) \,
\psi(t_0) \right) \, dxdydz  \nonumber  \\
&&\hskip-0.25in
=  \int_{t_0}^t
\int_{\Om}   (v, w) \cdot \left( (\phi_t, \psi_t) + \nu (\Dd \phi+\phi_{zz},
\Dd \psi+\psi_{zz})   \right)
\, dxdydz \; ds  \nonumber \\
&&\hskip-0.25in + \int_{t_0}^t \int_{\Om}  \left[  ( v \cdot
\nabla_h) (\phi, \psi)  \cdot (v, w) +w (\phi_z, \psi_z) \cdot (v,
w)  \right]  \, dxdydz +\int_{t_0}^t \int_{\Om}  \left[  ( f, g)
\cdot (v, w) \right]  \, dxdydz,    \label{WEAK}
\end{eqnarray*}
for every  $(\phi, \psi) \in \mathcal{V},$ and
almost every $t$, $t_0\in [0,T]$;

\item[(3)] the energy inequality:
\begin{eqnarray}
&&\hskip-.68in \frac{1}{2} \frac{d(\|v\|_2^2 +\|w\|_2^2)}{dt} +
\nu  \left[ \| \nabla_h v \|_2^2 +\|\nabla_h w \|_2^2
+\|v_z\|_2^2+\|w_z\|_2^2 \right]
 \leq  \int_{\Om} (f, g) \cdot (v, w)\; dxdydz.  \label{ENG}
\end{eqnarray}
\end{itemize}
Moreover, a weak solution is called strong solution of
(\ref{EQ-1})--(\ref{INIT-2}) on $[0,T]$ if, in addition, it
satisfies
\begin{eqnarray*}
&&(v, w)  \in C([0,T], V) \cap L^2([0,T], H^2 (\Om)).
\end{eqnarray*}

\begin{remark} \label{R-1}
Notice that one can extend $v$ and $w$
to be periodic, vertically, by, first, setting
\begin{eqnarray*}
&&\hskip-.8in
v(x,y,z, t) = v (x,y, -z, t),  \quad z \in (-1, 0) \\
&&\hskip-.8in
w(x,y,z,t) = - w (x,y,-z, t),  \quad z \in (-1, 0),
\end{eqnarray*}
and then, setting
\[
v(x,y,z+2, t) = v (x,y, z, t) \quad w(x,y,z+2,t) = w (x,y,z, t) \quad
\forall z \in \mathbb{R}.
\]
Similarly,  $p$, $f$, and $v_0$  can be extended to be even periodic functions
in the $z-$vertical with period 2,
 and $g$ and $w_0$ are odd periodic functions in the $z-$vertical with period 2.
In this way, the system (\ref{EQ-1})--(\ref{INIT-2}) can be treated
as the $3D$ Navier--Stokes equations with periodic boundary
condition.   In other words the system (\ref{EQ-1})--(\ref{INIT-2})
is a special case of the $3D$ Navier--Stokes equations on
$\mathbb{R}^3$ with periodic boundary condition. Then, by standard
procedure for the $3D$ Navier--Stokes equations with periodic
boundary condition (see, e.g., \cite{CF88}, \cite{DJ95},
\cite{LADY03}, \cite{RS93}, \cite{SO01A}, \cite{TT84}, \cite{TT95})
one can show that there exists, global in time, a Leray--Hopf weak
solution to the system (\ref{EQ-1})--(\ref{INIT-2}) if $(v_0, w_0)
\in H$. Furthermore, one can show the short time existence of the
strong solution if $(v_0, w_0) \in V$.

\end{remark}

\vskip0.1in

In this paper, we provide sufficient conditions on the pressure which guarantee the
global existence of the strong solution to
the 3D Navier--Stokes equations in infinite channel subject to the boundary conditions
(\ref{B-1})--(\ref{B-2}).  Several authors (see, for example, \cite{BG02}, \cite{CL01},
\cite{DV00}, \cite{SS02},  \cite{ZY04}, \cite{ZY05})
have studied the question of the global regularity of the $3D$ Navier--Stokes equation
by providing
sufficient conditions on the pressure.
Specifically, the authors of \cite{SS02} have shown that if $|u|^2 + 2p$ is bounded from below or from
above then the weak solution is strong solution.
The authors of \cite{CL01} have shown that if $p$ satisfies
\begin{equation}
  p\in L^{\aa}(0,T;L^{\beta}(\mathbb{R}^3))\quad\text{with}
\quad \frac{2}{\aa}+\frac{3}{\beta}<2,
\quad \text{for}\quad \beta>\frac{3}{2},    \label{CHE}
\end{equation}
then the weak solution is strong solution. This result has been improved
in \cite{BG02} by assuming
\begin{equation}
  p\in L^{\aa}(0,T;L^{\beta}(\mathbb{R}^3))\quad\text{with}\quad
\frac{2}{\aa}+\frac{3}{\beta}=2,
\quad \text{for}\quad \beta>\frac{3}{2}.   \label{GADI}
\end{equation}
In \cite{ZY05} the author has established the global regularity by assuming that
\begin{equation}
 \nabla_h p, \; p_z \in L^{\aa}(0,T;L^{\beta}(\mathbb{R}^3))\quad
\text{with}\quad \frac{2}{\aa} +\frac{3}{\beta}\leq 3,
\quad \text{for}\quad \aa>\frac{2}{3}, \;\;  \beta>1.   \label{ZHOU}
\end{equation}
Here, we will show the existence of the strong solutions of the
system (\ref{EQ-1})--(\ref{INIT-2}) on interval $[0, T]$ provided the vertical
derivative of the pressure satisfies
\begin{equation}
  p_z\in L^{\aa}([0,T]; L^{\beta}(\Om)) \quad \text{with} \quad \aa > 3, \quad
\text{and} \quad \beta > 2.    \label{CON}
\end{equation}
Let us observe that the quantity that appears in the Navier--Stokes
equations is $(\nabla_h p, p_z)$ and not $p$ itself. Therefore, the
conditions from \cite{ZY05} and our conditions seem to be natural.
Furthermore, it is  worth mentioning that the conditions (\ref{CON})
seem not to be comparable in any direct way with (\ref{CHE}),
(\ref{GADI}) and (\ref{ZHOU}). Nevertheless, our tools are totally
different from  those in \cite{BG02}, \cite{CL01} and \cite{ZY05}.
We prove our result by using the methods we established in
\cite{CT05} in which we proved the global well--posedness of the
$3D$ viscous primitive equations. Finally, notice that the
conditions (\ref{SP}) and (\ref{GADI}) are scaling invariant when
the equation is considered on whole $\mathbb{R}^3$. Namely,
\[
\left( \int_0^T \|u_{\la}\|_{\beta}^{\aa} \; dt\right)^{\frac{1}{\aa}}
= \left( \int_0^T \|u \|_{\beta}^{\aa} \; dt\right)^{\frac{1}{\aa}},
\qquad \mbox{when} \quad \frac{2}{\aa}+\frac{3}{\beta}=1,
\]
and
\[
 \left( \int_0^T \|p_{\la}\|_{\beta}^{\aa} \; dt\right)^{\frac{1}{\aa}}
= \left( \int_0^T \|p \|_{\beta}^{\aa} \; dt\right)^{\frac{1}{\aa}},
\qquad \mbox{when} \quad \frac{2}{\aa}+\frac{3}{\beta}=2,
\]
where $u_{\la} = \la u(\la x, \la y, \la z, \la^2 t)$ and $p_{\la} =
\la^2 u(\la x, \la y, \la z, \la^2 t)$ which is also a solution to
the Navier--Stokes equation if $(u, p)$ is a solution to the
Navier--Stokes equation in whole $\mathbb{R}^3$. However, our
condition is not scaling invariant. The reason may be because either
our results are not optimal, or because our condition involves  only
 one partial derivative of the pressure term and our result is
limited to channel flows and not in the whole space $\mathbb{R}^3$,
where the scaling argument is applied.

The plan of this paper is as follows. In section \ref{S-2}, we
reformulate the system (\ref{EQ-1})--(\ref{INIT-2}),  introduce our
notations, and recall some well--known useful inequalities. Section
\ref{S-3} is the main section, and it is devoted for the regularity
of solutions.

%%%%%%%%%%%%%%%%%%%%%%%%%%%%%%%%%%%%%%%%%%%%%%%%%%%%%%%%%%%%%%%%%%%%%%%%
%%%%%%%%%%%%%%%%%%%%  SECTION 2 %%%%%%%%%%%%%%%%%%%%%%%%%%%%%%%%%%%%%%%%
%%%%%%%%%%%%%%%%%%%%%%%%%%%%%%%%%%%%%%%%%%%%%%%%%%%%%%%%%%%%%%%%%%%%%%%%

\section{Preliminaries and Functional Setting}    \label{S-2}

In this section we introduce a  new equivalent formulation of
(\ref{EQ-1})--(\ref{INIT-2}). Following the ideas introduced in
\cite{CT05} we integrate the equation (\ref{EQ-3}) in the $z$
direction and by (\ref{B-1}) we obtain
\begin{equation}
w(x,y,z,t) =  - \int_{0}^z \nabla_h \cdot v(x,y, \xi,t) d\xi,
\label{DIV-1}
\end{equation}
and
\begin{equation}
\int_{0}^1 \nabla_h \cdot v(x,y, z,t) dz = \nabla_h \cdot
\int_{0}^1 v(x,y, z,t) dz =0. \label{DIV}
\end{equation}
For every function $\tt(x,y,z)$, we denote by
\begin{eqnarray}
&&\hskip-.8in \overline{\tt}(x,y)= \int_0^{1} \tt (x,y,z)\; dz.
\label{AVER}
\end{eqnarray}
and
\[
\widetilde{\tt}= \tt -\overline{\tt}.
\]
Following the geophysical terminology we denote  by
\begin{equation}
\overline{v} (x, y)= \frac{1}{h} \int_{-h}^0  v(x,y, \xi) d\xi,
\qquad \mbox{in } \; M, \label{V--B}
\end{equation}
 the barotropic mode.  We will also  denote by
\begin{equation}
\widetilde{v} = v - \overline{v},    \label{V--T}
\end{equation}
the baroclinic mode, that is the fluctuation about the barotropic
mode.
Notice that
\begin{eqnarray}
&&\hskip-.8in \overline{\widetilde{\tt}} =0. \label{ZERO}
\end{eqnarray}
By substituting (\ref{DIV-1}) into (\ref{EQ-1}), we reach
\begin{eqnarray}
&&\hskip-.68in \frac{\pp v}{\pp t} -\nu \Dd v - \nu v_{zz}+
(v\cdot \nabla_h) v - \left( \int_{0}^z \nabla_h \cdot v(x,y,
\xi,t) d\xi \right)  \frac{\pp v}{\pp  z}  +  \nabla_h p = f,
\label{VV}
\end{eqnarray}
and
\begin{eqnarray}
&&\hskip-.68in \nabla_h \cdot \overline{v} = 0.
\label{VDIV}
\end{eqnarray}

Furthermore, by taking the average of equations (\ref{VV}) in the $z$
direction, over the interval $(0, 1)$ and using the boundary
conditions (\ref{B-1}), we reach
\begin{eqnarray}
&&\hskip-.68in \frac{\pp \overline{v}}{\pp t} -\nu \Dd
\overline{v} + \overline{(v\cdot \nabla_h) v - \left( \int_{0}^z
\nabla_h \cdot v(x,y, \xi,t) d\xi \right) \frac{\pp v}{\pp   z}} +
\nabla_h \overline{p}  = \overline{f}.   \label{VA}
\end{eqnarray}
As a result of (\ref{DIV}), (\ref{ZERO}), and
integration by parts and using the boundary
conditions (\ref{B-1}),
the nonlinear term in (\ref{VA}) gives
\begin{eqnarray}
&&\hskip-.268in
 \overline{(v\cdot \nabla_h) v - \left( \int_{0}^z
\nabla_h \cdot v(x,y, \xi,t) d\xi \right) \frac{\pp v}{\pp   z}}=
(\overline{v} \cdot \nabla_h ) \overline{v} + \overline{
\left[(\widetilde{v} \cdot \nabla_h) \widetilde{v}
 + (\nabla_h \cdot \widetilde{v}) \; \widetilde{v}\right]}.     \label{I--1}
\end{eqnarray}
By subtracting  (\ref{VA}) from (\ref{VV}) and using (\ref{I--1})
we obtain
\begin{eqnarray}
&&\hskip-.68in \frac{\pp \widetilde{v}}{\pp t} -\nu\Dd
\widetilde{v}-\nu \widetilde{v}_{zz} + (\widetilde{v} \cdot
\nabla_h) \widetilde{v} - \left( \int_{0}^z \nabla_h \cdot
\widetilde{v}(x,y, \xi,t) d\xi \right) \frac{\pp
\widetilde{v}}{\pp  z}  \nonumber    \\
&&\hskip-.58in  +(\widetilde{v} \cdot \nabla_h ) \overline{v}+
(\overline{v} \cdot \nabla_h) \widetilde{v}  - \overline{
\left[(\widetilde{v} \cdot \nabla_h) \widetilde{v}  + (\nabla_h
\cdot \widetilde{v}) \; \widetilde{v}\right]} + \nabla_h
\widetilde{p} = \widetilde{f}.  \label{EQ4}
\end{eqnarray}
In addition, $\widetilde{v}$
satisfies the boundary conditions:
\begin{eqnarray}
&&\hskip-.68in
\widetilde{v}(x+1,y,z)=\widetilde{v}(x,y+1,z)=v(x,y,z),   \label{EQ6}  \\
&&\hskip-.68in \left. \frac{\pp \widetilde{v} }{\pp z}
\right|_{z=0} = 0, \; \left. \frac{\pp \widetilde{v} }{\pp z}
\right|_{z=1} = 0. \label{EQ66}
\end{eqnarray}

\vskip0.2in

For convenience, we recall the following Gagiliardo--Nirenberg, Sobolev and
Ladyzhenskaya inequalities in $\mathbb{R}^2$ (see, e.g., \cite{AR75},
\cite{CF88}, \cite{GA94}, \cite{LADY})
\begin{eqnarray}
&&\hskip-.68in \| \phi\|_{L^{\aa} (M)} \leq C_{\aa} \| \phi
\|_{L^2(M)}^{2/{\aa}} \| \phi \|_{H^1(M)}^{\frac{{\aa}-2}{{\aa}}},
\label{SI-2}
\end{eqnarray}
for every $\phi \in H^1(M), 2\leq \aa < \infty$,
and the following Gagiliardo--Nirenberg,
Sobolev and Ladyzhenskaya inequalities in $\mathbb{R}^3$
(see, e.g., \cite{AR75}, \cite{CF88}, \cite{GA94}, \cite{LADY})
\begin{eqnarray}
&&\hskip-.68in \| \psi \|_{L^{\aa}(\Om)} \leq C_{\aa} \| \psi
\|_{L^2(\Om)}^{\frac{6-\aa}{2\aa}} \| \psi
\|_{H^1(\Om)}^{\frac{3({\aa}-2)}{2{\aa}}},  \label{SI1}
\end{eqnarray}
for every $\psi\in H^1(\Om), 2\leq \aa \leq 6.$  Here $C_{\aa}$ denote constants
which are scale invariant.
Also, by (\ref{SI-2}) we get
\begin{eqnarray}
&&\hskip-.68in \| \phi \|_{L^{\beta}(M)} = \| |\phi|^{\aa/2}
\|_{L^{2\beta/\aa}(M)}^{2/\aa}
\leq C \| |\phi|^{\aa/2} \|_{L^2(M)}^{2/\beta} \| |\phi|^{\aa/2}
\|_{H^1(M)}^{\frac{2(\beta-\aa)}{\aa\,\beta}}  \nonumber  \\
&&\hskip-.68in \leq C  \| \phi \|_{L^\aa(M)}^{\aa/\beta} \; \left(
\int_M |\phi|^{\aa-2} \left| \nabla_h \phi \right|^2 \; dxdy
\right)^{\frac{(\beta-\aa)}{\aa\,\beta}} + \| \phi \|_{L^\aa(M)},
\label{TWE}
\end{eqnarray}
for every $\phi \in H^1(M),$ and $\beta > \aa.$  Also, we recall the integral
version of Minkowsky inequality for the $L^\beta$ spaces, $\beta\geq 1$.
Let $\Om_1 \subset \mathbb{R}^{m_1}$ and
 $\Om_2 \subset \mathbb{R}^{m_2}$ be two measurable sets, where
$m_1$ and $m_2$ are two positive integers. Suppose that
$f(\xi,\eta)$ is measurable over $\Om_1 \times \Om_2$. Then,
\begin{equation}
\hskip0.35in \left[ { \int_{\Om_1} \left( \int_{\Om_2} |f(\xi,\eta)| d\eta
\right)^{\beta} d\xi } \right]^{1/\beta}
\leq \int_{\Om_2} \left( \int_{\Om_1} |f(\xi,\eta)|^{\beta} d\xi
\right)^{1/\beta} d\eta.
\label{MKY}
\end{equation}

%%%%%%%%%%%%%%%%%%%%%%%%%%%%%%%%%%%%%%%%%%%%%%%%%%%
%%%%%%%%%%%%%%%%%%%%%%%%%%%%%%%%%%%%%%%%%%%%%%
%%%%%%%%%%%%%%%   SECTION 3 %%%%%%%%%%%%%%
%%%%%%%%%%%%%%%%%%%%%%%%%%%%%%%%%%%%%%%%%%%%%%%%%%%%
%%%%%%%%%%%%%%%%%%%%%%%%%%%%%%%%%%%%%%%%%%%%%%%%%%%
%%%%%%%%%%%%%%%%%%%%%%%%%%%%%%%%%%%%%%%%%%%%%%%%%%%%
%%%%%%%%%%%%%%%%%%%%%%%%%%%%%%%%%%%%%%%%%%%%%%%%%%%

\vskip0.2in

\section{Existence of the Strong Solution} \label{S-3}

In this section we will prove the global existence of the strong
solution to the system (\ref{EQ-1})--(\ref{INIT-2}) under the
assumption (\ref{CON}) on $p_z$.

\begin{theorem} \label{T-MAIN}
Suppose that $f, g \in H^1(\Om).$ For every $(v_0, w_0) \in
V,$ and  if $p_z \in L^{\aa}([0, T], L^{2q})$ with $T>0, \aa>3$ and $q>1$,
then there is a unique strong solution $((v, w), p)$ of the
system {\em (\ref{EQ-1})--(\ref{EQ-3})} on $[0, T]$.
\end{theorem}

\vskip0.05in

\begin{proof}
In Remark \ref{R-1} we described an argument that guarantees the
existence of a Leray--Hopf weak solution and short time existence of
the strong solutions. Suppose that $((v,w), p)$ is the strong
solution with initial value $(v_0, w_0) \in V$ such that $(v, w) \in
C([0, \mathcal{T}^*), V)  \cap L^2 ([0, \mathcal{T}^*), H^2(\Om)),$
where $[0, \mathcal{T}^*), \, \text{for}\, \mathcal{T}^* \leq T,$ is
the maximal interval of existence. If $\mathcal{T}^* = T$, then
there is nothing to prove. Next, we would like to show that certain
norms of this strong solution remain finite for all the time, up to
$\mathcal{T}^*$, provided the condition (\ref{CON}) is valid. In
this way we show that $\mathcal{T}^*$ is, at least, equal to $T$.
Namely, the strong solution $((v,w), p)$ exists on $[0, T]$. By the
energy inequality (\ref{ENG}) we have (see, for example,
\cite{CF88}, \cite{SO01A}, \cite{TT84} for details)
\begin{eqnarray}
&&\hskip-.68in \|v(t)\|_2^2 +\|w(t)\|_2^2 \leq K_{11},   \label{K11}
\end{eqnarray}
and
\begin{eqnarray}
&&\hskip-.68in \nu \int_0^t \left[ \| \nabla_h v (s)\|_2^2
+\|\nabla_h w(s)\|_2^2 +\|v_z(s)\|_2^2+\|w_z(s)\|_2^2 \right] \;
ds \leq K_{12}(t),  \label{K12}
\end{eqnarray}
where
\begin{eqnarray}
&&\hskip-.68in K_{11}= \frac{\|f\|_2^2+\|g\|_2^2}{\nu^2 \la_1^2} +
\|v_0\|_2^2+\|w_0\|_2^2. \label{K-11}   \\
&&\hskip-.68in K_{12}(t)= \frac{(\|f\|_2^2+\|g\|_2^2)\; t}{\nu \la_1} +
\|v_0\|_2^2+\|w_0\|_2^2. \label{K-12}
\end{eqnarray}

\subsection{{\text \bf $\|\widetilde{v}\|_r, 3 <  r <4$} estimates}

Taking the  inner product of the  equation (\ref{EQ4}) with
$|\widetilde{v}|^{r-2} \widetilde{v}$ in $L^2(\Om)$,  we get
\begin{eqnarray*}
&&\hskip-.168in \frac{1}{r} \frac{d \| \widetilde{v} \|_{r}^{r}
}{d t} + \nu \int_{\Om} \left(|\nabla_h \widetilde{v}|^2
|\widetilde{v}|^{r-2} + (r-2)\left|\; \nabla_h |\widetilde{v}|\;
\right|^2 |\widetilde{v}|^{r-2} + |\widetilde{v}_z|^2
|\widetilde{v}|^{r-2}
+(r-2) \left|\;\pp_z |\widetilde{v}| \;\right|^2 |\widetilde{v}|^{r-2} \right) \;
dxdydz    \\
&&\hskip-.165in = - \int_{\Om} \left(  (\widetilde{v} \cdot
\nabla_h) \widetilde{v} - \left( \int_{-h}^z \nabla_h \cdot
\widetilde{v}(x,y, \xi,t) d\xi \right)  \frac{\pp
\widetilde{v}}{\pp  z} +(\widetilde{v} \cdot \nabla_h )
\overline{v}+ (\overline{v} \cdot \nabla_h) \widetilde{v}  \right.  \\
&&\hskip-.01in
 \left.  -
\overline{ (\widetilde{v} \cdot \nabla_h) \widetilde{v} +
(\nabla_h \cdot \widetilde{v}) \; \widetilde{v}} +  \nabla_h
\widetilde{p} - \widetilde{f}\right) \cdot |\widetilde{v}|^{r-2}
\widetilde{v} \; dxdydz.
\end{eqnarray*}
By integration by parts and using the boundary conditions
(\ref{EQ6}) and (\ref{EQ66}) we get
\begin{eqnarray}
&&\hskip-.065in   - \int_{\Om} \left( (\widetilde{v} \cdot
\nabla_h) \widetilde{v} - \left( \int_{0}^z \nabla_h \cdot
\widetilde{v}(x,y, \xi,t) d\xi \right)  \frac{\pp
\widetilde{v}}{\pp  z} \right) \cdot |\widetilde{v}|^{r-2}
\widetilde{v} \; dxdydz =0. \label{D6-1}
\end{eqnarray}
Moreover, (\ref{VDIV}) and the boundary condition (\ref{EQ66})
give us
\begin{eqnarray}
&&\hskip-.065in \int_{\Om} (\overline{v} \cdot \nabla_h)
\widetilde{v}\cdot |\widetilde{v}|^{r-2} \widetilde{v} \;dxdydz =0.
\label{D6-3}
\end{eqnarray}
Thus, from (\ref{D6-1})--(\ref{D6-3}) we have
\begin{eqnarray*}
&&\hskip-.168in \frac{1}{r} \frac{d \| \widetilde{v} \|_{r}^{r}
}{d t} + \nu \int_{\Om} \left(|\nabla_h \widetilde{v}|^2
|\widetilde{v}|^{r-2} + (r-2)\left|\;\nabla_h |\widetilde{v}|\;
\right|^2 |\widetilde{v}|^{r-2} + |\widetilde{v}_z|^2
|\widetilde{v}|^{r-2}
+(r-2) \left|\;\pp_z |\widetilde{v}| \;\right|^2 |\widetilde{v}|^{r-2} \right) \;
dxdydz    \\
&&\hskip-.165in
 = - \int_{\Om} \left( (\widetilde{v} \cdot \nabla_h ) \overline{v}  - \overline{
(\widetilde{v} \cdot \nabla_h) \widetilde{v} + (\nabla_h \cdot
\widetilde{v}) \; \widetilde{v}}   + \nabla_h   \widetilde{p}-
\widetilde{f}  \right) \cdot |\widetilde{v}|^{r-2} \widetilde{v}
\; dxdydz.
\end{eqnarray*}
Notice that by integration by parts and using boundary conditions
(\ref{EQ6}) and (\ref{EQ66}) we have
\begin{eqnarray*}
&&\hskip-.165in
 - \int_{\Om} \left[
(\widetilde{v} \cdot \nabla_h ) \overline{v}  - \overline{
(\widetilde{v} \cdot \nabla_h) \widetilde{v} + (\nabla_h \cdot
\widetilde{v}) \; \widetilde{v}} + \nabla_h \widetilde{p}
\right] \cdot |\widetilde{v}|^{r-2} \widetilde{v} \; dxdydz  \\
&&\hskip-.165in
 = \int_{\Om} \left[
(\nabla_h \cdot \widetilde{v})  \; \overline{v} \cdot
|\widetilde{v}|^{r-2} \widetilde{v} + (\widetilde{v} \cdot
\nabla_h ) (|\widetilde{v}|^{r-2} \widetilde{v}) \cdot
\overline{v}
 - \sum_{k,j=1}^2 \left[ \overline{ \widetilde{v^k}   \widetilde{v^j}} \;  \pp_{x_k}
 (|\widetilde{v}|^{r-2} \widetilde{v^j}) \right] +
 \widetilde{p} \left(\nabla_h \cdot (|\widetilde{v}|^{r-2}
\widetilde{v}) \right) \right]  \; dxdydz.
\end{eqnarray*}
Observe that since $\overline{\widetilde{p}} =0$, we have the
Poincar\'{e} inequality
\begin{eqnarray}
&&\hskip-.065in
|\widetilde{p} | \leq \int_{-h}^0 |p_z| \; dz.
\label{PPP-3}
\end{eqnarray}
Therefore, from all the above and  by Cauchy--Schwarz and H\"{o}lder inequalities
we obtain
\begin{eqnarray*}
&&\hskip-.368in \frac{1}{r} \frac{d \| \widetilde{v} \|_{r}^{r}
}{d t} + \nu \int_{\Om} \left(|\nabla_h \widetilde{v}|^2
|\widetilde{v}|^{r-2} + (r-2)\left|\;\nabla_h |\widetilde{v}|\;
\right|^2 |\widetilde{v}|^{r-2} + |\widetilde{v}_z|^2
|\widetilde{v}|^{r-2}
+(r-2) \left|\;\pp_z |\widetilde{v}| \;\right|^2 |\widetilde{v}|^{r-2} \right) \;
dxdydz   \\
&&\hskip-.365in \leq C \int_{M} \left[ |\overline{v}| \int_{0}^1
|\nabla_h \widetilde{v}| \; | \widetilde{v}|^{r-1} \; dz \right]
\; dxdy  \\
&&\hskip-.26in +C \int_{M} \left[ \int_{0}^1 |\widetilde{v}|^2 \;
dz  \int_{0}^1 |\nabla_h \widetilde{v}| \;
| \widetilde{v}|^{r-2} \; dz \right] \; dxdy \\
&&\hskip-.265in +C \int_{M} \left[ \int_0^1 |p_z| \; dz \;
\int_{0}^1 |\nabla_h \widetilde{v}| \; | \widetilde{v}|^{r-2} \;
dz \right]
\; dxdy  + \|\widetilde{f}\|_r \; \|\widetilde{v}\|_r^{r-1}  \\
&&\hskip-.365in \leq C \int_{M} \left[ |\overline{v}| \left(
\int_{0}^1 |\nabla_h \widetilde{v}|^2 \; | \widetilde{v}|^{r-2} \;
dz \right)^{1/2} \left( \int_{0}^1  | \widetilde{v}|^r \; dz
\right)^{1/2} \right]
\; dxdy  \\
&&\hskip-.26in +C \int_{M} \left[ \int_{0}^1 |\widetilde{v}|^2 \;
dz \left( \int_{0}^1 |\nabla_h \widetilde{v}|^2 \; |
\widetilde{v}|^{r-2} \; dz \right)^{1/2}
\left( \int_{0}^1  | \widetilde{v}|^{r-2} \; dz \right)^{1/2} \right] \; dxdy \\
&&\hskip-.26in +C \int_{M} \left[ \int_0^1 |p_z| \; dz \; \left(
\int_{0}^1 |\nabla_h \widetilde{v}|^2 \; | \widetilde{v}|^{r-2} \;
dz \right)^{1/2} \left( \int_{0}^1  | \widetilde{v}|^{r-2} \; dz
\right)^{1/2} \right]
\; dxdy + \|\widetilde{f}\|_r \; \|\widetilde{v}\|_r^{r-1}  \\
&&\hskip-.365in \leq C \|\overline{v}\|_{L^4(M)}
 \left( \int_{\Om}  |\nabla_h \widetilde{v}|^2 \;
 | \widetilde{v}|^{r-2} \; dxdydz \right)^{1/2}
\left( \int_M \left( \int_{0}^1  | \widetilde{v}|^r \; dz
\right)^{2}
\; dxdy \right)^{1/4}   \\
&&\hskip-.26in +C \left( \int_M \left( \int_{0}^1  |
\widetilde{v}|^2 \; dz \right)^{\frac{r+2}{2}} \; dxdy
\right)^{\frac{2}{r+2}} \left( \int_{\Om}  |\nabla_h
\widetilde{v}|^2 \; | \widetilde{v}|^{r-2} \; dxdydz
\right)^{\frac{1}{2}} \left( \int_M \left( \int_{0}^1  |
\widetilde{v}|^{r-2} \; dz \right)^{\frac{r+2}{r-2}}
\; dxdy \right)^{\frac{r-2}{2(r+2)}} \\
&&\hskip-.26in +C \|\overline{|p_z|}\|_{L^{2q}(M)} \left(
\int_{\Om} |\nabla_h \widetilde{v}|^2 \; | \widetilde{v}|^{r-2} \;
dxdydz \right)^{1/2} \left( \int_M \left( \int_{0}^1  |
\widetilde{v}|^{r-2} \; dz \right)^{q^{\prime}} \; dxdy
\right)^{1/{2q^{\prime}}},
\end{eqnarray*}
where $1/q+1/q^{\prime} =1.$ By using Minkowsky  inequality (\ref{MKY}),  we
get
\begin{eqnarray*}
&&\hskip-.168in \left( \int_M \left( \int_{0}^1  |
\widetilde{v}|^r \; dz \right)^{2} \; dxdy \right)^{1/2} \leq  C
\int_{0}^1 \left( \int_M  | \widetilde{v}|^{2r} \; dxdy
\right)^{1/2}  \; dz.
\end{eqnarray*}
By virtue of (\ref{TWE}), we have
\begin{eqnarray*}
&&\hskip-.168in \int_M  | \widetilde{v}|^{2r}  \; dxdy \leq C_0
\int_M | \widetilde{v}|^{r}  \; dxdy \int_M | \widetilde{v}|^{r-2}
|\nabla_h \widetilde{v} |^2  \; dxdy + \left( \int_M |
\widetilde{v}|^{r}  \; dxdy \right)^2.
\end{eqnarray*}
Thus, by Cauchy--Schwarz inequality we obtain
\begin{eqnarray}
&&\hskip-.168in \left( \int_M \left( \int_{0}^1  |\widetilde{v}|^r \; dz \right)^{2}
\; dxdy \right)^{1/2} \leq  C
\| \widetilde{v}\|_r^{r/2}  \left( \int_{\Om} |\widetilde{v}|^{r-2}
|\nabla_h \widetilde{v} |^2 \; dxdydz \right)^{1/2} + \|
\widetilde{v}\|_r^{r}. \label{M1}
\end{eqnarray}
Similarly, by (\ref{MKY}) and (\ref{TWE}) and (\ref{SI-2}), we
also obtain
\begin{eqnarray}
&&\hskip-.168in \left( \int_M \left( \int_{0}^1  | \widetilde{v}|^2 \; dz
\right)^{(r+2)/2} \; dxdy \right)^{2/(r+2)}
\leq C \int_{0}^1 \left( \int_M  | \widetilde{v}|^{2+r}
\; dxdy \right)^{2/(r+2)}  \; dz   \nonumber  \\
&&\hskip-.168in \leq C \int_0^1 \left[ \|
\widetilde{v}\|_{L^r(M)}^{\frac{2r}{r+2}} \left( \int_{M}
|\widetilde{v}|^{r-2} |\nabla_h \widetilde{v} |^2 \; dxdy
\right)^{\frac{4}{r(r+2)}} + \| \widetilde{v}\|_{L^r(M)}^{2}
\right]
\; dz   \nonumber  \\
&&\hskip-.168in \leq C \int_0^1 \left[ \|
\widetilde{v}\|_{L^2(M)}^{\frac{2}{r}}\|\nabla_h
\widetilde{v}\|_{L^2(M)}^{\frac{r-2}{r}} \|
\widetilde{v}\|_{L^r(M)}^{\frac{r-2}{r+2}} \left( \int_{M}
|\widetilde{v}|^{r-2} |\nabla_h \widetilde{v} |^2 \; dxdy
\right)^{\frac{4}{r(r+2)}} + \| \widetilde{v}\|_{L^r}^{2} \right]
\; dz \nonumber  \\
&&\hskip-.168in \leq C  \|
\widetilde{v}\|_2^{\frac{2}{r}}\|\nabla_h
\widetilde{v}\|_2^{\frac{r-2}{r}} \|
\widetilde{v}\|_r^{\frac{r-2}{r+2}} \left( \int_{\Om}
|\widetilde{v}|^{r-2} |\nabla_h \widetilde{v} |^2 \; dxdydz
\right)^{\frac{4}{r(r+2)}} + \| \widetilde{v}\|_r^{2}, \label{M2}
\end{eqnarray}
and
\begin{eqnarray}
&&\hskip-.168in \left( \int_M \left( \int_{0}^1  |
\widetilde{v}|^{r-2} \; dz \right)^{(r+2)/(r-2)} \; dxdy
\right)^{\frac{r-2}{r+2}} \leq C \int_{0}^1 \left( \int_M  |
\widetilde{v}|^{2+r}
\; dxdy \right)^{(r-2)/(r+2)}  \; dz   \nonumber  \\
&&\hskip-.168in \leq C \int_0^1 \left[ \|
\widetilde{v}\|_{L^r(M)}^{\frac{r(r-2)}{r+2}} \left( \int_{M}
|\widetilde{v}|^{r-2} |\nabla_h \widetilde{v} |^2 \; dxdy
\right)^{\frac{2(r-2)}{r(r+2)}} + \|
\widetilde{v}\|_{L^r(M)}^{r-2} \right]
\; dz   \nonumber  \\
&&\hskip-.168in  \leq C  \| \widetilde{v}\|_r^{\frac{r(r-2)}{r+2}}
\left( \int_{\Om} |\widetilde{v}|^{r-2} |\nabla_h \widetilde{v} |^2
\; dxdydz \right)^{\frac{2(r-2)}{r(r+2)}} + \|
\widetilde{v}\|_r^{r-2}, \label{M22}
\end{eqnarray}
Thanks to (\ref{SI-2}) and (\ref{MKY}), we conclude
\begin{eqnarray}
&&\hskip-.168in \left( \int_M \left( \int_{0}^1  |
\widetilde{v}|^{r-2} \; dz \right)^{\frac{r}{r-2}} \; dxdy \right)^{\frac{r-2}{r}}
\leq  C \int_{0}^1 \left( \int_M  | \widetilde{v}|^{r}
\; dxdy \right)^{\frac{r-2}{r}}  \; dz    \nonumber  \\
&&\hskip-.168in \leq  C  \int_{0}^1  \left[ \|
\widetilde{v}\|_{L^2(M)}^{\frac{2(r-2)}{r}} \;
 \|\nabla_h \widetilde{v}\|_{L^2(M)}^{\frac{(r-2)^2}{r}}
 + \| \widetilde{v}\|_{L^2(M)}^{r-2} \right] \; dz
\leq  C   \|
\widetilde{v}\|_2^{\frac{2(r-2)}{r}} \;
 \|\nabla_h \widetilde{v}\|_2^{\frac{(r-2)^2}{r}}   + \| \widetilde{v}\|_{L^2(M)}^{r-2}.
\label{M3}
\end{eqnarray}
Moreover,
\begin{eqnarray}
&&\hskip-.168in \left( \int_M \left( \int_{0}^1  |
\widetilde{v}|^{r-2} \; dz \right)^{q^{\prime}} \; dxdy \right)^{1/q^{\prime}} \leq  C
\int_{0}^1 \left( \int_M  | \widetilde{v}|^{(r-2)q^{\prime}}
\; dxdy \right)^{1/q^{\prime}}  \; dz    \nonumber  \\
&&\hskip-.168in \leq  C  \int_{0}^1  \left[ \|
\widetilde{v}\|_{L^r(M)}^{r/q^{\prime}} \left( \int_{M}
|\widetilde{v}|^{r-2} |\nabla_h \widetilde{v} |^2 \; dxdy
\right)^{\frac{(r-2)q^{\prime}-r}{rq^{\prime}}}   + \|
\widetilde{v}\|_{L^r(M)}^{(r-2)} \right] \; dz
\nonumber  \\
&&\hskip-.168in \leq C   \| \widetilde{v}\|_r^{r/q^{\prime}}
\left( \int_{\Om} |\widetilde{v}|^{r-2} |\nabla_h \widetilde{v} |^2
\; dxdydz \right)^{\frac{(r-2)q^{\prime}-r}{rq^{\prime}}}   + \|
\widetilde{v}\|_r^{(r-2)}. \label{M4}
\end{eqnarray}
Therefore, (\ref{M1})--(\ref{M4}) and (\ref{SI-2}) give
\begin{eqnarray*}
&&\hskip-.168in \frac{d \| \widetilde{v} \|_{r}^{r} }{d t} + \nu
\int_{\Om} \left(|\nabla_h \widetilde{v}|^2 |\widetilde{v}|^{r-2}
+ \left|\nabla_h |\widetilde{v}| \right|^2 |\widetilde{v}|^{r-2} +
|\widetilde{v}_z|^2 |\widetilde{v}|^{r-2}
+ \left|\pp_z |\widetilde{v}| \right|^2 |\widetilde{v}|^{r-2} \right) \; dxdydz    \\
&&\hskip-.165in \leq C \|\overline{v}\|_2^{1/2} \; \|\nabla_h
\overline{v}\|_2^{1/2} \| \widetilde{v}\|_r^{r/2}
 \left( \int_{\Om}  |\nabla_h \widetilde{v}|^2 \; | \widetilde{v}|^4 \; dxdydz
 \right)^{3/4}
+C \|\overline{v}\|_2^{1/2} \; \|\nabla_h \overline{v}\|_2^{1/2}
\| \widetilde{v}\|_r^{r} + \|f\|_r\; \|\widetilde{v}\|_r^{r-1} + C
\|\widetilde{v}\|_r^r
   \\
&&\hskip-.065in +C \| \widetilde{v}\|_2^{\frac{2}{r}} \|\nabla_h
\widetilde{v}\|_2^{\frac{r-2}{r}} \|
\widetilde{v}\|_r^{\frac{r-2}{2}} \left( \int_{\Om}
|\widetilde{v}|^{r-2} |\nabla_h \widetilde{v} |^2 \; dxdydz
\right)^{\frac{r+2}{2r}}  \\
&&\hskip-.065in +C \|p_z\|_{2q} \;  \|
\widetilde{v}\|_r^{\frac{r}{2q^{\prime}}} \left( \int_{\Om}
|\nabla_h \widetilde{v}|^2 \; | \widetilde{v}|^{r-2} \; dxdydz
\right)^{1-\frac{r+2q^{\prime}}{2rq^{\prime}}} + C \|p_z\|_{2q}\;
\| \widetilde{v}\|_r^{(r-2)/2} \; \left( \int_{\Om} |\nabla_h
\widetilde{v}|^2 \; | \widetilde{v}|^{r-2} \; dxdydz
\right)^{1/2}.
\end{eqnarray*}
By Young's and Cauchy--Schwarz inequalities we have
\begin{eqnarray*}
&&\hskip-.168in
  \frac{d \| \widetilde{v} \|_{r}^{r} }{d t} + \nu
\int_{\Om} \left(|\nabla_h \widetilde{v}|^2 |\widetilde{v}|^{r-2}
+ \left|\nabla_h |\widetilde{v}| \right|^2 |\widetilde{v}|^{r-2} +
|\widetilde{v}_z|^2 |\widetilde{v}|^{r-2}
+ \left| \pp_z |\widetilde{v}| \right|^2 |\widetilde{v}|^{r-2} \right) \; dxdydz   \\
&&\hskip-.165in \leq C\left[\; 1+ \|\overline{v}\|_2^{2} \;
\|\nabla_h \overline{v}\|_2^{2}  +  \|
\widetilde{v}\|_2^{\frac{4}{r-2}} \|\nabla_h \widetilde{v}\|_2^2
\right] \| \widetilde{v}\|_r^{r} +C \|p_z\|_{2q}^{r}+C\|f\|_r^r.
\end{eqnarray*}
Thanks to  Gronwall inequality, we get
\begin{eqnarray}
&&\hskip-.68in \| \widetilde{v} (t)\|^r_r + \nu \int_0^t
\int_{\Om} \left[ |\nabla_h \widetilde{v}|^2
|\widetilde{v}|^{r-2}+ |\widetilde{v}_z|^2 |\widetilde{v}|^{r-2}
\right] \; dxdydz \leq K_{R},   \label{K-R}
\end{eqnarray}
where
\begin{eqnarray}
&&\hskip-.68in K_R = e^{C T+K_{11} K_{12}(T) + K_{11}^{\frac{2}{r-2}}K_{12}(T)}
\left[\; 1+
\|v_0\|_{H^1(\Om)}^6 + \int_0^T \|p_z (s)\|_{2q}^{r}\;
ds+\|f\|_r^r T \right].   \label{K6}
\end{eqnarray}
It is worth mentioning that the above estimate is also valid for
$2 \leq r < 4.$ However, one need $r>3$ in order to get the
following $H^1$ estimate.

\subsection{$H^1$ estimates}

Before we show the global $H^1$ bound, let us
prove the following Lemma.

\begin{lemma} \label{LLL}
Let $\phi \in H^1(\Om), \psi \in L^2(\Om),$ and $v, r$ be as in
Theorem \ref{T-MAIN}. Then,
\begin{eqnarray}
&&\hskip-.165in
\int_{\Om} |v| \, |\phi| \, |\psi| \; dxdydz \leq  \ee
\left(\|\nabla \phi\|_2^2
+\|\phi_z\|_2^2 + \|\psi\|_2^2 \right) \\
&& + C_{\ee}\left[ \|\widetilde{v}\|_r^{\frac{2r}{r-3}}
 + \|\widetilde{v}\|_r^2  +  \left(
1+ \|\overline{v}\|_2^2 \right) \; \left( \|\overline{v}\|_2^2+
\|\nabla_h \overline{v}\|_2^2 \right)\right]  \|\phi\|_2^2    \label{LL}
\end{eqnarray}
for every $\ee >0.$
\end{lemma}

\begin{proof}
Notice that
\begin{eqnarray*}
&&\hskip-.65in
 \int_{\Om} |v| \, |\phi| \, |\psi| \; dxdydz \leq
 \int_{\Om} \left( |\widetilde{v}| +|\overline{v}| \right) |\phi| \, |\psi| \;
 dxdydz.
\end{eqnarray*}
By H\"{o}lder inequality and (\ref{SI-2}), we obtain
\begin{eqnarray*}
&&\hskip-.65in
 \int_{\Om} |\widetilde{v}| \, |\phi| \, |\psi| \;
 dxdydz \leq \|\widetilde{v}\|_r \, \|\phi\|_{\frac{2r}{r-2}} \, \|\psi\|_2 \\
&&\hskip-.65in
 \leq C \|\widetilde{v}\|_r \, \|\phi\|_2^{\frac{r-3}{r}} \, \left( \|\nabla \phi\|_2
+\|\phi_z\|_2 \right)^{\frac{3}{r}} \|\psi\|_2
 + \|\widetilde{v}\|_r  \, \|\phi\|_2 \, \|\psi\|_2.
\end{eqnarray*}
By Young's inequality, we reach
\begin{eqnarray*}
&&\hskip-.65in
 \int_{\Om} |\widetilde{v}| \, |\phi| \, |\psi| \;
 dxdydz \leq  \frac{\ee}{2} \left(\|\nabla \phi\|_2^2
+\|\phi_z\|_2^2 + \|\psi\|_2^2 \right)+ C_{\ee}  \left(
\|\widetilde{v}\|_r^{\frac{2r}{r-3}}
 + \|\widetilde{v}\|_r^2  \right)\; \|\phi\|_2^2.
\end{eqnarray*}
On  the other hand, by applying the same method for proving Proposition 2.2 in
\cite{CT03}, we get
\begin{eqnarray*}
&&\hskip-.65in
 \int_{\Om} |\overline{v}| \, |\phi| \, |\psi| \;
 dxdydz \leq C \left(\|\overline{v}\|_2^{1/2}\|\nabla_h
 \overline{v}\|_2^{1/2}+\|\overline{v}\|_2\right)
 \, \left(\|\phi\|_2^{1/2} \,  \|\nabla_h \phi\|_2
+\|\phi\|_2\right) \; \|\psi\|_2.
\end{eqnarray*}
Again, by Young's inequality, we obtain
\[
\int_{\Om} |v| \, |\phi| \, |\psi| \; dxdydz \leq  \frac{\ee}{2}
\left(\|\nabla_h \phi\|_2^2 +\|\psi\|_2^2 \right)+C_{\ee} \left(
1+ \|\overline{v}\|_2^2 \right) \; \left( \|\overline{v}\|_2^2+
\|\nabla_h \overline{v}\|_2^2 \right) \|\phi\|_2^2.
\]
Therefore, (\ref{LL}) holds.

\end{proof}

Next, let show the $H^1$ norm of the strong solution $(v, w)$ is bounded.
Taking the  inner product of the  equation (\ref{EQ-1}) with $-\Dd v -v_{zz}$
and the  equation (\ref{EQ-2}) with $-\Dd w -w_{zz}$ in
$L^2$,  and using the fact that the Stokes operator is same as the Laplacian operator under
periodic boundary conditions, we obtain
\begin{eqnarray*}
&&\hskip-.22in \frac{1}{2} \frac{d \left(\|\nabla_h v \|_2^2+\|
v_z \|_2^2 +\|\nabla_h w \|_2^2+\| w_z \|_2^2 \right) }{d t} + \nu
\left( \|\Dd v\|_2^2+ 2 \|\nabla_h v_z\|_2^2 + \|v_{zz}\|_2^2 +
\|\Dd w\|_2^2+ 2 \|\nabla_h
w_z\|_2^2 + \|w_{zz}\|_2^2    \right) \\
&&\hskip-.2in = \int_{\Om} \left[  (v \cdot \nabla_h) v +wv_z -f
\right] \cdot \left(\Dd v+ v_{zz}  \right) \; dxdydz+  \int_{\Om}
\left[  v \cdot \nabla_h w +ww_z -g \right] \left(\Dd w+ w_{zz}
\right) \; dxdydz.
\end{eqnarray*}
By integration by parts and using the boundary conditions
(\ref{B-1}) and (\ref{B-2}), we obtain
\begin{eqnarray*}
&&\hskip-.065in    \int_{\Om} wv_z \cdot \left(\Dd v+ v_{zz} \right) \;
dxdydz   \\
&&\hskip-.065in = -  \int_{\Om} \left( ((\nabla_h w \cdot
\nabla_h v) \cdot v_z + w_z   |v_z|^2 \right)  \; dxdydz  \\
&&\hskip-.065in  =   \int_{\Om} \left( \nabla_h w_z  \cdot
\nabla_h v +\nabla_h w \cdot  \nabla_h v_z + w_{zz} \; v_z +w_z \; v_{zz}
\right) \cdot v \; dxdydz,
\end{eqnarray*}
and also
\begin{eqnarray*}
&&\hskip-.065in    \int_{\Om} ww_z \left(\Dd w+ w_{zz} \right) \;
dxdydz   \\
&&\hskip-.065in = -  \int_{\Om} \left( |\nabla_h w|^2 w_z + w
\nabla_h w_z \cdot \nabla_h w+ \frac{1}{2} w_z^3\right)
\; dxdydz  \\
&&\hskip-.065in  = - \frac{1}{2} \int_{\Om} \left(|\nabla_h w|^2
+w_z^2\right) w_z \; dxdydz
\\
&&\hskip-.065in  = \frac{1}{2} \int_{\Om} \left(|\nabla_h w|^2
+w_z^2\right)  \left( \nabla_h \cdot v \right) \; dxdydz \\
&&\hskip-.065in = - \frac{1}{2} \int_{\Om} \nabla_h
\left(|\nabla_h w|^2 +w_z^2\right) \cdot v \; dxdydz.
\end{eqnarray*}
Then, from the above, we get
\begin{eqnarray*}
&&\hskip-.22in \frac{1}{2} \frac{d \left(\|\nabla_h v \|_2^2+\|
v_z \|_2^2 +\|\nabla_h w \|_2^2+\| w_z \|_2^2 \right) }{d t} + \nu
\left( \|\Dd v\|_2^2+ 2 \|\nabla_h v_z\|_2^2 + \|v_{zz}\|_2^2 +
\|\Dd w\|_2^2+ 2 \|\nabla_h
w_z\|_2^2 + \|w_{zz}\|_2^2    \right)  \\
&&\hskip-.2in = \int_{\Om} \left[  (v \cdot \nabla_h) v  -f
\right] \cdot \left(\Dd v+ v_{zz}  \right) \; dxdydz + \int_{\Om}
\left( \nabla_h w_z  \cdot
\nabla_h v +\nabla_h w \cdot  \nabla_h v_z + w_{zz} \; v_z +w_z \; v_{zz}
\right) \cdot v \; dxdydz  \\
&&\hskip-.2in +  \int_{\Om} \left[  v \cdot \nabla_h w  -g
\right] \left(\Dd w+ w_{zz} \right) \; dxdydz - \frac{1}{2}
\int_{\Om} \nabla_h \left(|\nabla_h w|^2 +w_z^2\right) \cdot
v \; dxdydz.
\end{eqnarray*}
By applying Lemma \ref{LLL} with some small enough $\ee$ and the
Cauchy--Schwarz inequality, we obtain
\begin{eqnarray*}
&&\hskip-.268in  \frac{d \left(\|\nabla_h v \|_2^2+\| v_z \|_2^2
+\|\nabla_h w \|_2^2+\| w_z \|_2^2 \right) }{d t} + \nu \left(
\|\Dd v\|_2^2+ 2 \|\nabla_h v_z\|_2^2 + \|v_{zz}\|_2^2 + \|\Dd
w\|_2^2+ 2 \|\nabla_h
w_z\|_2^2 + \|w_{zz}\|_2^2    \right)  \\
&&\hskip-.165in \leq C \left( 1+
\|\widetilde{v}\|_r^{\frac{2r}{r-3}} + \|\overline{v}\|_2^4 \;
+ \|\nabla_h \overline{v}\|_2^4 \right) \left(\|\nabla_h v
\|_2^2+\| v_z \|_2^2 +\|\nabla_h w \|_2^2+\| w_z \|_2^2 \right) +
\|f\|^2_2 +\|g\|^2_2.
\end{eqnarray*}
Thanks to Gronwall inequality, we obtain
\begin{eqnarray*}
&&\hskip-.68in \|\nabla_h v (t)\|_2^2+\| v_z (t)\|_2^2 +\|\nabla_h
w(t) \|_2^2+\| w_z (t)\|_2^2 \\
&&\hskip-.6in +\nu \int_0^t \left( \|\Dd v(s)\|_2^2+ 2 \|\nabla_h
v_z(s)\|_2^2 + \|v_{zz}(s)\|_2^2 + \|\Dd w(s)\|_2^2+ 2 \|\nabla_h
w_z(s)\|_2^2 + \|w_{zz}(s)\|_2^2  \right) \; ds \leq K_2,
\end{eqnarray*}
where
\begin{eqnarray}
&&\hskip-.68in K_2 = e^{CT+K_{11}\;(T+ K_{12}(T))  +K_R^{2/(r-3)} \; T} \left[
\|v_0\|_{H^1(\Om)} +\|w_0\|_{H^1(\Om)}+ \|f\|^2_2 +\|g\|^2_2
\right].   \label{K2}
\end{eqnarray}
Therefore, the $H^1$ norm of the solution remains bounded on the
maximal interval of existence $[0,\mathcal{T}^*)$. This completes
the proof of theorem.
\end{proof}

\noindent
\section*{Acknowledgements}
E.S.T. would like to thank the   Bernoulli Center of the \'{E}cole
Polytechnique F\'{e}d\'{e}ral de Lausanne  where part of this work
was completed. This work was supported in part by the NSF grant No.
DMS-0504619, the ISF grant No. 120/60 and the BSF grant No. 2004271.

\end{document}